\documentclass[12pt, a4paper,twoside]{article}
\usepackage{a4,amssymb,amsfonts,verbatim,pstricks,pst-plot,xypic,times}
\usepackage{yfonts}

\usepackage[english]{babel}
\usepackage{amssymb,latexsym}
\usepackage{amssymb}
\usepackage{latexsym}

\newcommand{\SO}{\mathrm{SO}}
\newcommand{\M}{\mathbb{M}}

\newcommand{\paa}{{\partial}}
\newcommand{\HH}{{\mathbb H}}
\newcommand{\CC}{{\mathbb C}}
\newcommand{\DD}{{\mathbb D}}
\newcommand{\Pe}{\mathrm{P}}
\newcommand{\RR}{{\mathbb R}}
\renewcommand{\Re}{\mathrm{Re}}
\newcommand{\ZZ}{{\mathbb Z}}

\newcommand{\SSS}{{\mathbb S}}

\newcommand{\Ekt}{\mathbb{E}(\kappa,\tau)}
\newcommand{\NN}{{\mathbb N}}

\newcommand{\cZ}{{\mathcal Z}}

\renewcommand{\phi}{\varphi}

\newcommand{\half}{\frac{1}{2}}

\newcommand{\cercle}{\mathbb{S}}
\newcommand{\grad}{\mathrm{grad}}

\newcommand{\Spin}{\mathrm{Spin}}
\renewcommand{\Re}{\mathrm{Re}}
\newcommand{\End}{\mathrm{End}}

\newcommand{\tr}{\mathrm{tr}}
\newcommand{\<}{\left\langle}       
\renewcommand{\>}{\right\rangle}

\newcommand{\Spinc}{\mathrm{Spin^c}}

\newcommand{\DDM}{\widetilde D}

\newcommand{\Chi}{\mathfrak{X}}

\newcommand{\id}{\mathrm{Id}} 
\newcommand{\vol}{\mathrm{vol}}

\newtheorem{example}{Examples}[section]
\newtheorem{thm}{Theorem}[section]
\newtheorem{lemma}[thm]{Lemma}
\newtheorem{prop}[thm]{Proposition}

\newtheorem{remark}[thm]{Remark}
\newtheorem{remarks}[thm]{Remarks}
\newtheorem{definition}[thm]{Definition}
\newtheorem{notation}[thm]{Notation}
\newtheorem{exabout:ample}[thm]{Example}

\begin{document}
\title{The Spin$^c$ Dirac Operator on Hypersurfaces and Applications}
\author{{\small\bf Roger NAKAD}\footnote{Max Planck Institute for Mathematics, Vivatsgasse 7, 53111 Bonn, Germany. E-mail: nakad@mpim-bonn.mpg.de}\ \ and
{\small\bf{ Julien ROTH}}\footnote{Laboratoire d'Analyse et de Math\'ematiques appliqu\'ees, Universit\'e Paris-Est Marne-la-Vall\'ee, 77454 Marne-la-Vall\'ee. E-mail: julien.roth@univ-mlv.fr}}
\maketitle
\begin{center}
{\bf Abstract}
\end{center} 

We extend to the eigenvalues of the hypersurface $\Spinc$ Dirac operator well known lower and upper  bounds. Examples of limiting cases 
are then given. Futhermore, we prove
 a correspondence between the existence of a $\Spinc$ Killing spinor on homogeneous $3$-dimensional 
manifolds $\mathbb E^*(\kappa, \tau)$ with $4$-dimensional isometry group and 
isometric immersions of $\mathbb E^*(\kappa, \tau)$
into the complex space form $\M^4(c)$ of constant holomorphic sectional curvature $4c$, for some $c\in \RR^*$. 
As applications, we show the non-existence  
of totally umbilic surfaces in $\mathbb E^*(\kappa, \tau)$ and we give necessary and sufficient geometric 
conditions to immerse a $3$-dimensional Sasaki manifold into $\M^4(c)$.
\section{Introduction}
It is well known that the spectrum of the Dirac operator on hypersurfaces of a $\Spin$ manifold detects informations 
on the geometry of such manifolds and their hypersurfaces (\cite{Am1, Am2, bar, HMZ01, HMZ02, HMR02}). For example, 
O. Hijazi, S. Montiel and X. Zhang \cite{HMZ01} proved that on the compact boundary  $M^n$ of a Riemannian compact $\Spin$ manifold 
$\cZ^{n+1}$ of dimension $n+1$ and with nonnegative scalar curvature, the first positive eigenvalue $\lambda_1$ of the Dirac operator satisfies
\begin{eqnarray}
\lambda_1 \geq \frac n2 \inf_M H,
\label{1lo}
\end{eqnarray}
where $H$ denotes the mean curvature of $M$, assumed to be nonnegative. Equality holds if and only if $H$ is constant and
 every eigenspinor associated with $\lambda_1$ is the restriction to $M$ of a parallel spinor on $\cZ$ (and so $\cZ$ is 
Ricci-flat). As application of the limiting case, they gave an elementary proof of the famous Alexandrov theorem
 \cite{HMZ01}: {\it the only compact embedded hypersurface in $\RR^{n+1}$ of constant mean curvature is the sphere $\mathbb{S}^n$ of dimension $n$}.\\\\
Assume now that $M^n$ is a closed hypersurface of $\cZ^{n+1}$. Evaluating the Rayleigh quotient applied to a parallel or Killing
 spinor field coming from $\cZ$, C. B\"ar \cite{bar}  derived an upper bound for the eigenvalues of the Dirac operator on 
$M$ by using the min-max principle. More precisely, there are at least $\mu$ eigenvalues $\lambda_1, \cdots, \lambda_\mu$ 
 of the Dirac operator on $M$ satisfying
\begin{eqnarray}
\lambda_j^2 \le n^2\alpha^2 + \frac{n^2}{4 \ \vol(M)}\int_M H^2 dv,
\label{inegualitee}
\end{eqnarray}
where $\vol(M)$ is the volume of $M$, $dv$ is the volume form of the manifold $M$, $\alpha$ is the Killing number
 ($\alpha =0$ if the ambient spinor field is parallel) and $\mu$ is the dimension of the space of parallel or Killing spinors. 
\\\\
Recently, $\Spinc$ geometry  became a field of active research with the advent
of Seiberg-Witten theory \cite{KM94, Wit94, Sei-Wit94}.  Applications of the Seiberg-Witten theory to $4$-dimensional geometry and topology
are already notorious (\cite{Don96, LeB95, LeB96, GL98}). From an intrinsic point of view, $\Spin$, almost complex, 
complex, K\"{a}hler, Sasaki
and some classes of CR manifolds have a canonical $\Spinc$ structure. The complex projective space $\CC P^m$ is always
 $\Spinc$ but not $\Spin$ if $m$ is even. Nowadays, and from the extrinsic point of
vue, it seems that it is more natural to work with $\Spinc$ structures rather than
Spin structures. Indeed, O. Hijazi, S. Montiel and F. Urbano \cite{HMU1} constructed on
K\"{a}hler-Einstein manifolds with positive scalar curvature, $\Spinc$ structures carrying
K\"{a}hlerian Killing spinors. The restriction of these spinors to minimal Lagrangian
submanifolds provides topological and geometric restrictions on these submanifolds.
In \cite{NR, nakadthese}, and via $\Spinc$ spinors, the authors gave an elementary proof for a {\it Lawson type correspondence} between constant 
mean curvature surfaces of $3$-dimensional homogeneous manifolds with $4$-dimensional isometry group. 
We point out that, using $\Spin$ spinors, we cannot prove this {\it Lawson type correspondence}. Moreover, they characterized
isometric immersions of a $3$-dimensional almost contact metric manifold $M$ into the complex space form
 by the existence  of a $\Spinc$ structure on $M$ carrying a special spinor called a generalized Killing spinor.\\\\
In the first part of this paper and using the $\Spinc$ Reilly inequality, we extend the lower bound (\ref{1lo}) to the first positive 
eigenvalue of the Dirac operator defined on the compact boundary of a $\Spinc$ manifold. Examples of the limiting case 
are then given where the equality case cannot occur if we consider the $\Spin$ Dirac operator on these examples. Also, by restriction of parallel and Killing $\Spinc$ spinors, we extend the upper bound
 (\ref{inegualitee}) to the eigenvalues of the Dirac operator defined on  a closed hypersurface of $\Spinc$ manifolds. 
Examples of the limiting  case are also given. \\\\
In the second part, we study  $\Spinc$ structures on $3$-dimensional homogeneous manifolds $\mathbb{E}^*(\kappa, \tau)$ with 
$4$-dimensional isometry group. It is well known that the only  complete simply connected $\Spinc$ manifolds admitting 
real Killing spinor other than the $\Spin$ manifolds are the non-Einstein Sasakian manifolds endowed with their canonical
 or anti-canonical $\Spinc$ structure \cite{19}. Since $\mathbb{E}^*(\kappa, \tau)$ are non-Einstein Sasakian manifolds 
\cite{ein}, the canonical and the anti-canonical $\Spinc$ structure carry real Killing spinors. In \cite{NR}, the authors
 proved that this canonical (resp. this anti-canonical) $\Spinc$ structure on $\mathbb{E}^*(\kappa, \tau)$ is the lift of the canonical (resp. the
 anti-canonical) $\Spinc$ structure on $\M^2(\kappa)$ via the submersion $\mathbb{E}^*(\kappa, \tau) \longrightarrow \M^2(\kappa)$, where $\M^2(\kappa)$ denotes the simply connected $2$-dimensional manifold with 
constant curvature $\kappa$ . Moreover, 
they proved that the Killing constant of the real  Killing spinor field is equal to $\frac{\tau}{2}$. Here, we reprove the 
existence of a Killing spinor for the canonical and the anti-canonical $\Spinc$ structure. This proof is based on the existence of an isometric embedding of $\mathbb{E}^*(\kappa, \tau)$ into the complex projective space or 
the complex hyperbolic space (see Proposition \ref{im-inv}). Conversely, from the existence of a Killing spinor on $\mathbb{E}^*(\kappa, \tau)$, 
we prove the existence of an isometric  immersion of $\mathbb{E}^*(\kappa, \tau)$ into the complex space form $\M^4(c)$ of 
constant holomorphic sectional curvature $4c$, for some $c \in \RR^*$ (see Proposition \ref{im}). Since every non-Einstein Sasaki manifold has a  $\Spinc$ structure with a Killing spinor, it is natural to ask if this last 
result remains true for any $3$-dimensional Sasaki manifold. Indeed, every  simply connected
 non-Einstein Sasaki manifold can be immersed 
into $\M^4(c)$ for some $c\in \RR^*$, providing that the scalar curvature is constant (see Theorem \ref{sasaki-immersion}). Finally, we make use of the existence of a 
 Killing spinor on $\mathbb{E}^*(\kappa, \tau)$ to calculate some eigenvalues of Berger spheres  endowed 
with differents $\Spinc$ structures. By restriction of this Killing spinor to any surface of $\mathbb{E}^*(\kappa, \tau)$, 
we give a $\Spinc$ proof for the 
non-existence of totally umbilic surfaces in $\mathbb{E}^*(\kappa, \tau)$ (see Theorem \ref{tum}) proved already by R. Souam and E. Toubiana \cite{ST}. 
\section{Preliminaries}
In this section, we briefly introduce basic notions concerning the Dirac operator on $\Spinc$ manifolds (with or without 
boundary) and their hypersurfaces. Details can be found in \cite{6}, \cite{mon}, \cite{16}, \cite{hijazi4} and \cite{bar}.\\\\
{\bf The Dirac operator on Spin$^c$ manifolds.} We consider an oriented Riemannian manifold  $(M^n, g)$  of dimension $n$ with or without boundary and denote by $\SO M$ the 
$\SO_n$-principal bundle over $M$ of positively oriented orthonormal frames. A $\Spinc$ structure of $M$ is is given by an $\cercle^1$-principal bundle $(\cercle^1 M ,\pi,M)$ of some Hermitian line bundle $L$ and a $\Spin_n^c$-principal bundle $(\Spinc M,\pi,M)$ which is a $2$-fold covering of the $\SO_n\times\cercle^1$-principal bundle $\SO M\times_{M}\cercle^1 M$ compatible with the group covering
$$0 \longrightarrow \ZZ_2 \longrightarrow \Spin_n^c = \Spin_n\times_{\ZZ_2}\cercle^1 \longrightarrow \SO_n\times\cercle^1 \longrightarrow 0.$$
The bundle $L$ is called the auxiliary line bundle associated with the $\Spinc$ structure. If $A: T(\cercle^1 M)\longrightarrow i\RR$ is 
a connection 1-form on $\cercle^1 M$, its (imaginary-valued) curvature will be denoted by $F_A$, whereas we shall define a real 
$2$-form $\Omega$ on $\cercle^1 M$ by $F_A= i\Omega$. We know
that $\Omega$ can be viewed as a real valued 2-form on $M$ \cite{6, koba1}. In this case, $i\Omega$ is the curvature form of the auxiliary line bundle $L$ \cite{6, koba1}.\\\\
Let $\Sigma M := \Spinc M \times_{\rho_n} \Sigma_n$ be the associated spinor bundle where $\Sigma_n = \CC^{2^{[\frac n2]}}$ and $\rho_n : \Spin_n^c
\longrightarrow  \End(\Sigma_{n})$ the complex spinor representation \cite{6, 16, nakadthese}. 
A section of $\Sigma M$ will be called a spinor field. This complex vector bundle is naturally endowed
 with a Clifford multiplication, denoted by ``$\cdot$'', $\cdot : \CC l(TM) \longrightarrow \End(\Sigma M)$
 which is a fiber preserving algebra morphism and with a natural Hermitian scalar product $< . , .>$ 
compatible with this Clifford multiplication \cite{mon, 6, hijazi4}. If $n$ is even, $\Sigma M = \Sigma^+ M \oplus \Sigma^- M$ can be
 decomposed into positive and negative spinors by the action of the complex volume element 
\cite{6, mon, hijazi4, nakadthese}. If such data are given, one can canonically define a covariant derivative 
$\nabla$ on $\Sigma M$ given, for all $X \in \Gamma(TM)$, by \cite{6, 16, hijazi4, nakadthese}:
\begin{eqnarray}\label{nnabla}
 \nabla_X \psi = X(\psi) + \frac 14 \sum_{j=1}^n e_j\cdot \nabla_Xe_j\cdot\psi + \frac i2 A(s_*(X))\psi,
\end{eqnarray}
where the second $\nabla$ is the Levi-Civita connection on $M$, $\psi = [\widetilde{b\times s}, \sigma]$ is a locally defined spinor field, $b=(e_1, \cdots, e_n)$ is a local oriented orthonormal tangent frame, 
$s : U \longrightarrow \cercle^1M$ is a local section of $\cercle^1M$,  $\widetilde{b\times s}$ is the lift of the local section $b\times s :
 U \longrightarrow \SO M \times_M \cercle^1 M$ to the $2$-fold covering and
  $X(\psi) = [\widetilde{b\times s}, X(\sigma)]$. For any other connection $A^{'}$ on $\cercle^1 M$, 
there exists a real $1$-form $\alpha$ on $M$ such that  $ A^{'} = A  + i\alpha$ \cite{6}. If we endow the $\SSS^1$-principal fiber bundle $\SSS^1 M$ with the connection $A^{'}$, there exists on $\Sigma M$ a covariant derivative $\nabla^{'}$ given by 
\begin{eqnarray}\label{nnnabla}
\nabla_X^{'}\psi = \nabla_X \psi + \frac{i}{2} \alpha(X)\psi,
\end{eqnarray}
for all $X\in \Gamma(TM)$ and $\psi \in \Gamma(\Sigma M)$. Moreover, the curvature $2$-form of $A^{'}$ is given by $F_{A^{'}} = F_{A} + id\alpha$. 
But $F_{A}$ 
(resp. $F_{A^{'}}$) can be viewed as an imaginary $2$-form on $M$ denoted by $i\Omega$ (resp. $i\Omega^{'}$). Thus, 
$i\Omega$ (resp. $i\Omega^{'})$ is the curvature of the auxiliary line bundle associated with the
 $\SSS^1$-principal fiber bundle $\SSS^1 M$ endowed with the connection $A$ (resp. $A^{'}$) and we have $i\Omega^{'} = i\Omega + id\alpha$. \\\\
The Dirac operator, acting on $\Gamma(\Sigma M)$, is a first order elliptic operator locally given by  $D =\sum_{j=1}^n e_j \cdot \nabla_{e_j},$
where $\{e_j\}_{j=1, \cdots, n}$ is any orthonormal local basis tangent to $M$. An important tool when examining the Dirac operator on $\Spinc$ manifolds is the Schr\"{o}dinger-Lichnerowicz formula \cite{6, 16}:
\begin{eqnarray}
D^2 = \nabla^*\nabla + \frac 14 S\; \id_{\Gamma (\Sigma M)}+ \frac{i}{2}\Omega\cdot,
\label{sl}
\end{eqnarray}
where $S$ is the scalar curvature of $M$, $\nabla^*$ is the adjoint of $\nabla$ with respect to the $L^2$-scalar product and $\Omega\cdot$ is the extension of the Clifford multiplication to differential forms. 
The Ricci identity is given, for all $X\in \Gamma(TM)$,  by
\begin{eqnarray}\label{rici}
\sum_{j=1}^n e_j\cdot \mathcal R(e_j, X)\psi = \frac 12 \mathrm{Ric} (X)\cdot\psi - \frac i2 (X\lrcorner \Omega)\cdot\psi,
\end{eqnarray}
for any spinor field $\psi$. Here $\mathrm{Ric}$ (resp. $\mathcal R$) denotes the Ricci tensor 
of $M$ (resp. the $\Spinc$ curvature associated with the connection $\nabla$) and $\lrcorner$ the interior product.
\\\\
A $\Spin$ structure can be seen as a $\Spinc$ structure with trivial auxiliary line bundle
 $L$ endowed with the trivial connection. Every almost complex manifold $(M^{2m}, g, J)$ 
of complex dimension $m$ has a canonical $\Spinc$ structure. In fact, the complexified cotangent bundle 
$T^*M\otimes \CC = \Lambda^{1,0} M \oplus \Lambda^{0,1}M$ 
decomposes into the $\pm i$-eigenbundles of the complex linear extension of the complex structure. Thus, the spinor bundle of the canonical $\Spinc$ structure is given by $$\Sigma M = \Lambda^{0,*} M =\oplus_{r=0}^m \Lambda^{0,r}M,$$
where $\Lambda^{0,r}M = \Lambda^r(\Lambda^{0,1}M)$ is the bundle of $r$-forms of type $(0, 1)$. The auxiliary line bundle of this canonical $\Spinc$ structure is given by  $L = (K_M)^{-1}= \Lambda^m (\Lambda^{0,1}M)$, where $K_M$ is the canonical bundle of $M$ \cite{6, mon, nakadthese}. Let $\ltimes$ be the K\"{a}hler form defined by the complex structure $J$, i.e. $\ltimes (X, Y)= g(X, JY)$ for all vector fields $X,Y\in \Gamma(TM).$ The auxiliary line bundle $L= (K_M)^{-1}$ has a canonical holomorphic connection induced from the Levi-Civita connection whose curvature form is given by $i\Omega = i\rho$, where $\rho$ is the Ricci $2$-form given by $\rho(X, Y) = \mathrm{Ric} (X, JY)$. 
For any other $\Spinc$ structure the spinorial bundle can be written as \cite{6, HMU1}:
$$\Sigma M = \Lambda^{0,*}M\otimes\mathcal L,$$
where $\mathcal L^2 = K_M\otimes L$ and $L$  is the auxiliary bundle associated with this $\Spinc$
structure. In this case, the $2$-form $\ltimes$ can be considered as an endomorphism of $\Sigma M$ via
 Clifford multiplication and 
we have the well-known orthogonal splitting $\Sigma M = \oplus_{r=0}^{m}\Sigma_rM,$
where $\Sigma_rM$ denotes the eigensubbundle corresponding 
to the eigenvalue $i(m-2r)$ of $\ltimes$, with complex rank $\Big(^m_k\Big)$. The bundle $\Sigma_r M$ correspond 
to $\Lambda^{0, r}M\otimes\mathcal L$. For the canonical $\Spinc$ structure, the subbundle $\Sigma_0M$ is trivial. 
Hence and when $M$ is a  K\"{a}hler manifold, this $\Spinc$ structure admits parallel spinors (constant functions) 
lying in $\Sigma_0M$ \cite{19}. Of course, we can define another $\Spinc$ structure for which the spinor bundle is 
given by 
$\Lambda^{*, 0} M =\oplus_{r=0}^m \Lambda^r (T_{1, 0}^* M)$ and the auxiliary line bundle by $K_M$. 
This $\Spinc$ structure will be called the anti-canonical $\Spinc$ structure.\\\\
{\bf Spin$^c$ hypersurfaces and the Gauss formula.}
Let $(M^n, g)$ be an $n$-dimensional oriented hypersurface isometrically immersed in a Riemannian $\Spinc$ manifold 
$(\cZ^{n+1}, g_{\cZ})$. The hypersurface $M$ inherts a $\Spinc$ structure from that on $\cZ$, and we have \cite{mon, bar, nakadthese, 2ana}:
$$\left\{
\begin{array}{l}
\Sigma \cZ_{|_M} \ \ \ \ \simeq \Sigma M\ \text{\ \ \ if\ $n$ is even,} \\
\Sigma^+ \cZ_{|_M}\ \ \simeq\Sigma M \ \text{\ \ \ if\ $n$ is odd.}
\end{array}
\right.
$$
Moreover Clifford multiplication by a vector field $X$, tangent to $M$, is given by 
$$X\bullet\phi = (X\cdot \nu\cdot \psi)_{|_M},$$
where $\psi \in \Gamma(\Sigma \cZ)$ (or $\psi \in \Gamma(\Sigma^+ \cZ)$ if $n$ is odd), $\phi$ is the restriction of $\psi$
 to $M$, ``$\cdot$'' is the Clifford multiplication on $\cZ$, ``$\bullet$'' that on $M$ and $\nu$ is the unit normal vector. 
When $n$ is odd, we can also get $\Sigma^-\cZ_{\vert_{M}} \simeq \Sigma M$. In this case, the Clifford multiplication by a 
vector field $X$ tangent to $M$ 
is given by $X\bullet\phi = - (X\cdot\nu\cdot\psi)_{\vert_M}$ and we have $\Sigma\cZ_{|_M} \simeq \Sigma M\oplus\Sigma M$. The connection 1-form defined on the restricted $\SSS^1$-principal bundle $(\SSS^1 M =: {\SSS^1 \cZ}_{|_M},\pi,M)$, is given by
$$A = {A^\cZ}_{|_M} : T(\SSS^1 M) = T({\SSS^1 \cZ})_{|_M} \longrightarrow i\RR.$$
Then the curvature 2-form $i\Omega$ on the $\SSS^1$-principal bundle $\SSS^1 M$ is given by $i\Omega = {i\Omega^\cZ}_{|_M}$, 
which can be viewed as an imaginary 2-form on $M$ and hence as the curvature form of the line bundle $L$, the restriction
 of the line bundle $L^\cZ$ to $M$. We denote by $\nabla^{\cZ}$ the spinorial Levi-Civita connection on $\Sigma \cZ$ and 
by $\nabla$ that on $\Sigma M$. For all $X\in \Gamma(TM)$ and for every
 spinor field $\psi \in \Gamma(\Sigma \cZ)$ (or $\psi \in \Gamma(\Sigma^+ \cZ)$ if $n$ is odd), we consider $\phi= \psi_{|_M}$ and  we get the following $\Spinc$ Gauss formula \cite{mon, bar, 2ana}:
\begin{equation}
(\nabla^{\cZ}_X\psi)_{|_M} =  \nabla_X \phi + \half II(X)\bullet\phi,
\label{spingauss}
\end{equation}
where $II$ denotes the Weingarten map with respect to $\nu$. Moreover, Let $D^\cZ$ and $D$ be the Dirac operators on $\cZ$ and $M$, after denoting by the same symbol any spinor and it's restriction to $M$, we have
\begin{eqnarray}
\DDM\phi = \frac{n}{2}H\phi -\nu\cdot D^\cZ\phi-\nabla^{\cZ}_{\nu}\phi ,
\label{diracgauss}
\end{eqnarray}
\begin{eqnarray}
\DDM(\nu\cdot\varphi) = -\nu\cdot\widetilde D \varphi,
\label{d1}
\end{eqnarray}
where $H = \frac1n \tr(II)$ denotes the mean curvature and 
$\DDM = D$ if $n$ is even and $\DDM=D \oplus(-D)$ if $n$ is odd.
\\\\{\bf Homogeneous 3-dimensional manifolds with 4-dimensional isometry group.}
We denote a $3$-dimensional homogeneous manifold with $4$-dimensional isometry group by $\Ekt$, $\kappa-4\tau^2 \neq 0$. 
It is a Riemannian fibration over a simply connected $2$-dimensional manifold $\M^2(\kappa)$ with 
constant curvature $\kappa$ and such that the fibers are geodesic. We denote by $\tau$ the bundle curvature,
 which measures the default of the fibration to be a Riemannian product. Precisely, we denote by $\xi$
 a unit vertical vector field, that is tangent to the fibers. If $\tau \neq 0$, the vector field $\xi$ is a 
Killing field and satisfies for all  vector field $X$, 
$$\nabla_X\xi=\tau X\wedge\xi,$$ 
where $\nabla$ is the Levi-Civita connection and $\wedge$ is the exterior product. In this case $\Ekt$ is denoted by $\mathbb{E}^*(\kappa, \tau)$.
When $\tau$ vanishes, we get a product manifold $\M^2(\kappa)\times\RR$. If $\tau\neq0$,
 these manifolds are of three types: they have the isometry group of the Berger spheres if $\kappa>0$, 
of the Heisenberg group $\mathrm{Nil}_3$ if $\kappa=0$ or of $\widetilde{\mathrm{PSL}_2(\RR)}$ if $\kappa<0$. Note that if $\tau =0$, then $\xi =\frac{\partial}{\partial t}$ is the unit vector field giving 
the orientation of $\RR$ in the product $\M^2 (\kappa) \times \RR$. 
The manifold $\mathbb{E}^*(\kappa, \tau)$  admits a local direct orthonormal frame $\{e_1, e_2, e_3\}$ with  $e_3 = \xi,$ and such that the Christoffel symbols $ \Gamma^k_{ij} = \< \nabla_{e_i}e_j, e_k\>$ are given by
\begin{eqnarray}\label{christoffel}
\left\lbrace  
\begin{array}{l}
{\Gamma}_{12}^3={\Gamma}_{23}^1=-{\Gamma}_{21}^3=-{\Gamma}_{13}^2=\tau,\\ \\
{\Gamma}_{32}^1=-{\Gamma}_{31}^2=\tau-\frac{\kappa}{2\tau}, \\ \\
{\Gamma}_{ii}^i={\Gamma}_{ij}^i={\Gamma}_{ji}^i={\Gamma}_{ii}^j=0,\quad\forall\,i,j\in\{1,2,3\}.
\end{array}
\right. 
\end{eqnarray}
We call $\{e_1, e_2, e_3=\xi\}$ the canonical frame of $\mathbb{E}^*(\kappa, \tau)$. Except the Berger spheres and with $\RR^3$, $\HH^3$,
$\SSS^3$ and the solvable group $\mathrm{Sol}_3$, the manifolds $\Ekt$ define the geometry of Thurston. The authors 
\cite{NR} proved that there exists on $\mathbb{E}^*(\kappa, \tau)$ a $\Spinc$ structure (the canonical $\Spinc$ structure) carrying a Killing spinor field 
$\psi$ of Killing constant $\frac{\tau}{2}$, i.e., a spinor field $\psi$ satisfying
$$\nabla_X\psi = \frac{\tau}{2} X\cdot\psi,$$
for all $X\in \Gamma(T\mathbb{E}^*(\kappa, \tau))$. Moreover, $\xi\cdot\psi = -i\psi$ and the curvature of the auxiliary line bundle is given by
\begin{eqnarray}
 i\Omega(e_1, e_2) = -i(\kappa-4\tau^2)\ \ \ \text{and}\ \ \ \ \ i\Omega(e_k, e_j)=0,
\end{eqnarray}
elsewhere in the canonical frame $\{e_1, e_2, \xi\}$. There exists also another $\Spinc$ structure (the anti-canonical $\Spinc$ structure) carrying a Killing spinor 
 field $\psi$ of Killing constant $\frac{\tau}{2}$ such that $\xi\cdot\psi = i\psi$ and the curvature of the auxiliary line bundle is given by
\begin{eqnarray}
 i\Omega(e_1, e_2) = i(\kappa-4\tau^2)\ \ \ \text{and}\ \ \ \ \ i\Omega(e_k, e_j)=0,
\end{eqnarray}
elsewhere in the canonical frame $\{e_1, e_2, \xi\}$.
\section{Lower and upper bounds for the eigenvalues of the  hypersurface Dirac operator}
We will extend the lower bound (\ref{1lo}) and the upper bound (\ref{inegualitee}) 
to the eigenvalues of the  hypersurface $\Spinc$ Dirac operator $\widetilde D$. Examples of the limiting cases are then given.
\subsection{Lower  bounds for the eigenvalues of the  hypersurface Dirac operator }
We assume that the manifold $\cZ^{n+1}$ is a $\Spinc$ manifold having a compact domain $\DD$ with compact boundary
 $M= \partial \DD$. Using suitable boundary conditions for the Dirac operator $D^\cZ$, we extend the lower
 bound (\ref{1lo}) to the first positive eigenvalue of the extrinsic hypersurface Dirac operator 
$\widetilde D$ on $M$ endowed with the induced $\Spinc$ structure.\\\\
Since $M$ is compact, the Dirac operator $\widetilde D$ has a discrete spectrum and we denote by $\pi_+: \Gamma(\Sigma M) \longrightarrow \Gamma(\Sigma M)$ the projection onto the subspace of 
$\Gamma(\Sigma M)$ spanned by eigenspinors corresponding to the nonnegative eigenvalues of $\widetilde D$.  This projection provides an
 Atiyah-Patodi-Singer type boundary conditions for the Dirac operator $D^\cZ$ of the domain $\DD$. It has been proved that this is a global self-adjoint elliptic condition \cite{HMZ2, HMZ01}.\\\\
It is not difficult to extend the $\Spin$ Reilly inequality (see \cite{HMZ2}, \cite{HMZ01}, \cite{HMZ02}, \cite{HMR02}) to $\Spinc$ manifolds. Indeed, for all spinor fields $\psi\in \Gamma(\Sigma\DD)$, we have
\begin{eqnarray}\label{ir}
\int_{\paa \DD}\Big(< \DDM\varphi,\varphi>-\frac{n}{2} H|\varphi|^2\Big)ds\geq\int_{\DD}\Big(\frac{1}{4} S^\cZ \vert\psi\vert^2 &+&<\frac i2 \Omega^\cZ\cdot\psi, \psi> \nonumber\\ &-&\frac{n}{n+1}|D^\cZ\psi|^2\Big)dv,
\end{eqnarray}
where $dv$ (resp. $ds$) is the Riemannian volume form of $\DD$ (resp. $\paa\DD $). Moreover
equality occurs if and only if the spinor field $\psi$ is a twistor-spinor, i.e., if and only if it satisfies $\Pe^{\cZ}\psi=0$, where $\Pe^{\cZ}$ is the twistor operator
acting on $\Sigma\cZ$ locally given, for all $X\in\Gamma(T\cZ)$, by $\Pe^{\cZ}_X\psi=\nabla_X^{\cZ}\psi+\frac{1}{n+1}X\cdot D^
\cZ \psi$. 
Now, we can state the main theorem of this section:
\begin{thm}
Let $(\cZ^{n+1}, g_\cZ)$ be a Riemannian $\Spinc$ manifold such that the operator  $S^\cZ+2i\Omega^\cZ\cdot$ is nonnegative. We consider $M^n$ a compact hypersurface with 
nonnegative mean curvature $H$ and bounding a compact domain $\DD$ in $\cZ$. Then, the first positive eigenvalue
 $\lambda_1$ of $\widetilde D$  satisfies
\begin{eqnarray}
\lambda_1 \geqslant \frac{n}{2}\inf_M H.
\label{lower}
\end{eqnarray}
Equality holds if and only if $H$ is constant and the eigenspace corresponding to $\lambda_1$ consists of the restrictions 
to $M$ of parallel spinors on the domain $\DD$.
 \label{thth}
\end{thm}
{\bf Proof.} Let $\varphi$ be an eigenspinor on $M$ corresponding to the first positive eigenvalue $\lambda_1 > 0$ of $\widetilde D$, i.e., $\widetilde D \varphi = \lambda_1  \varphi$ and $\pi_+ \varphi = \varphi$. The following boundary problem has a unique solution (see \cite{HMZ2}, \cite{HMZ01}, \cite{HMZ02} and \cite{HMR02})
$$\left\{
\begin{array}{l}
D^\cZ \psi = 0\ \ \ \ \ \ \ \text{on} \ \ \DD \\ 
\pi_+ \psi = \pi_+ \varphi= \varphi\ \ \ \ \text{on}\ \ M=\partial\DD.
\end{array}
\right.
$$
From the Reilly inequality (\ref{ir}), we get 
$$\int_M (\lambda_1 - \frac n2 H)\vert\psi\vert^2 ds \geq \int_\DD (\frac 14 S^\cZ\vert\psi\vert^2  +\frac{i}{2}<\Omega^\cZ\cdot\psi, \psi>) dv \geq 0,$$
which implies (\ref{lower}). If the equality case holds in (\ref{lower}), then $\psi$ is a harmonic spinor and a twistor spinor, hence parallel. Since $\pi_+ \psi = \varphi$ along the boundary, $\psi$ is a non-trivial parallel spinor and $\lambda_1 = \frac n2 H$. Futhermore, since $\psi$ is parallel, we deduce by  (\ref{diracgauss}) that $\widetilde D \varphi = \frac n2 H \varphi$. Hence we have $\varphi = \pi_+ \psi = \psi$. Conversely if $H$ is constant, the fact that the restriction to $M$ of a parallel spinor on $\DD$ is an eigenspinor with eigenvalue $\frac n2 H$ is a direct consequence of (\ref{diracgauss}).
\begin{example}
\label{Ex1}
A complete simply connected Riemannian $\Spinc$ manifold $\cZ^{n+1}$ carrying a parallel spinor
 field is isometric to the Riemannian product of a simply connected K\"{a}hler manifold $\cZ_1^{n_1}$ of complex
 dimension $m_1$ ($n_1 =2m_1$) and a simply connected $\Spin$ manifold $\cZ_2^{n_2}$ of dimension $n_2$ 
($n +1 = n_1 + n_2$) carrying a parallel spinor and the $\Spinc$ structure of $\cZ$ is the product of the canonical
 $\Spinc$ structure of $\cZ_1$ and the $\Spin$ structure of $\cZ_2$ \cite{19}. Moreover, if we assume that $\cZ_1$ is Einstein, then 
\begin{eqnarray}\label{conditionspinc}
i\Omega^{\cZ} (X, Y) = i\rho^{\cZ_1}(X_1, Y_1) = i\mathrm{Ric}^{\cZ_1}(X_1, JY_1) = i \frac{S^{\cZ^1}}{n_1}\ltimes(X_1, Y_1),
\end{eqnarray}
for every $X = X_1 + X_2, Y = Y_1 + Y_2 \in \Gamma(T\cZ)$ and where $J$ denotes the complex structure 
on $\cZ_1$. Moreover, if the Einstein manifold $\cZ_1$ is of positive scalar curvature, we have, 
for any spinor field $\psi \in \Gamma(\Sigma\cZ)$, 
\begin{eqnarray*}
 S^\cZ \vert\psi\vert^2 + 2i<\Omega^\cZ\cdot\psi, \psi> &=& S^{\cZ_1} \vert\psi\vert^2  +\frac{i}{m_1}  S^{\cZ_1} <\ltimes\cdot\psi, \psi>\\&=&
 S^{\cZ_1} \sum_{r=0}^{m_1} (1 - \frac{m_1-2r}{m_1}) \vert\psi_r\vert^2 =  S^{\cZ_1} \sum_{r=0}^{m_1} \frac{2r}{m_1} 
\vert\psi_r\vert^2 \geq 0.
\end{eqnarray*}
Finally, the first positive eigenvalue of the Dirac operator $\widetilde D$ of any compact hypersurface with nonnegative
 constant mean curvature $H$ and  bounding a compact domain $\DD$ in $\cZ = \cZ_1 \times\cZ_2$ satisfies the equality case 
in (\ref{lower}) for the restricted $\Spinc$ structure. Next, we will give some explicit examples. The
 Alexandrov theorem for $\SSS_+^2 \times \RR$ says that the only embedded compact 
surface with constant mean curvature $H > 0$ in  $\cZ =\cZ_1\times \cZ_2 = \SSS_+^2 \times \RR$ is the standard rotational 
sphere described in \cite{AR1, AR2, benoit}. Hence, the first positive eigenvalue of the Dirac operator $\widetilde D$ on the rotational sphere satisfies the equality 
case in (\ref{lower}). 
We consider the complex projective space $\CC P^m$ ($\cZ_2 = \{\emptyset\}$) endowed 
with the Einstein Fubini-Study metric and the canonical $\Spinc$ structure. 
The first positive eigenvalue of the Dirac operator $\widetilde D$ of any compact hypersurface $M$ with nonnegative
 constant mean curvature $H$ and  bounding a compact domain $\DD$ in $\CC P^m$ satisfies the equality case in (\ref{lower}). 
Compact embedded hypersurfaces  in $\CC P^m$ are examples of manifolds viewed as a boundary of some enclosed domain 
in $\CC P^m$.  As an example, we know that there exists an isometric embedding of $\mathbb E^*(\kappa, \tau)$ into $\M^4(\frac{\kappa}{4}-\tau^2)$ 
of constant mean curvature $H = \frac{\kappa-16\tau^{2}}{12\tau}$ \cite{CMC}. 
Here $\M^4(\frac{\kappa}{4}-\tau^2)$ denotes the complex space form of constant holomorphic sectional 
curvature $\kappa-4\tau^2$. We choose $\kappa > 16 \tau^2$ and $\tau > 0$, then $H$ is positive.
 In this case, $\mathbb E^*(\kappa, \tau)$ are Berger spheres (compact) and $\M^4$ is the complex projective space $\CC P^2$ 
of constant holomorphic sectional curvature $\kappa-4\tau^2 > 0$. The canonical $\Spinc$ structure on $\M^4$ 
carries a parallel spinor and hence the equality case in
 (\ref{lower}) is satisfied for the first positive eigenvalue of $\widetilde D$ defined on Berger spheres. Finally, we recall that $\SSS_+^2 \times \RR$ and $\CC P^m$ (when $m$ is odd), have also a
 unique
 $\Spin$ structure. Hence, Inequality (\ref{lower}) holds for the first positive eigenvalue of
 the $\Spin$ Dirac operator $\widetilde D$ defined on the rotational sphere or on any compact embedded hypersurface in 
$\CC P^m$ (when $m$ is odd). But, equality cannot
 occur since this unique $\Spin$ structure on $\SSS_+^2 \times \RR$ and on $\CC P^m$ does not carry a parallel spinor.
\end{example}
\subsection{Upper bounds for the eigenvalues of the Dirac operator}
A spinor field $\psi$ on a Riemannian $\Spinc$ manifold $\cZ^{n+1}$ is called a real Killing spinor with Killing constant 
$\alpha\in \RR$ 
if 
\begin{eqnarray}
 \nabla^{\cZ}_X \psi = \alpha \ X\cdot\psi,
\label{killing}
\end{eqnarray}
for all $X \in \Gamma(T\cZ)$. When $\alpha=0$, the spinor field $\psi$ is a parallel spinor. We define
$$\mu = \mu (\cZ, \alpha) := \dim_\CC\{ \psi, \psi \ \ \text{is a Killing spinor on}\ \  \cZ \ \ \text{with Killing constant}\ \  \alpha \}$$
\begin{thm}
Let $M$ be an $n$-dimensional closed oriented hypersurface isometrically immersed in a Riemannian $\Spinc$ manifold $\cZ$. 
We endow $M$ with the induced $\Spinc$ structure. For any $\alpha \in \RR$, there are at least $\mu (\cZ, \alpha)$ eigenvalues $\lambda_1,\ldots,
\lambda_\mu$ of the Dirac operator $\widetilde D$ on $M$ satisfying
\begin{eqnarray}
\label{upperr}
 \lambda_j^2 \le n^2\alpha^2 + \frac{n^2}{4\vol(M)}\int_M H^2 dv,
\end{eqnarray}
\label{princ}
where $H$ denotes the mean curvature of $M$.  If equality holds, then $H$ is constant.
\end{thm}
{\bf Proof.} First, note that the set of Killing spinors with Killing constant $\alpha$ is a vector space. Moreover, linearly 
independent Killing spinors are linearly
independent at every point, the space of restrictions of
Killing spinors on $\cZ$ to $M$, i.e.,
$$
\{\psi|_M,  \ \psi \mbox{ is a spinor on $\cZ$ satisfying }
\nabla_X^{\cZ}\psi = \alpha\  X\cdot\psi,
\hspace{0.4cm} \forall X\in \Gamma(T\cZ) \}
$$
is also $\mu$-dimensional. Now let $\psi$ be a Killing spinor on $\cZ$ with Killing constant
$\alpha\in\RR$. Killing spinors have constant length so we can assume that $|\psi| \equiv 1$. By definition, we have 
$D^\cZ \psi = -(n+1)\alpha\  \psi$, and hence using (\ref{diracgauss}) 
we get $\DDM \phi = n\alpha\ \nu\cdot\phi + \frac {n}{2} H \phi.$ We denote by $(., . ) = \Re \int_M <., ,> $ the real part of the $L^2$-scalar product. 
Now, we compute the Rayleigh quotient of $\DDM^2$
$$\frac{(\DDM^2\psi,\psi)}{(\psi,\psi)} =
\frac{(\DDM\psi,\DDM\psi)}{\vol(M)}  = \frac{(n\alpha\  \nu\cdot\psi + \frac {n}{2} H \psi,n\alpha\ \nu\cdot\psi + \frac {n}{2} H \psi)}{\vol(M)} 
= n^2 \alpha^2 + \frac{n^2}{4} \frac{\int_M H^2}{\vol (M)},$$
i.e., the Rayleigh quotient of $\DDM^2$ is bounded by
$n^2\alpha^2 + \frac{n^2}{4}\frac{\int_M H^2}{\vol(M)}$ on a
$\mu$-dimensional space of spinors on $M$. Hence, the Min-Max
principle implies the assertion. If equality holds, then the restriction to $M$ of every Killing spinor $\psi$ of Killing constant $\alpha$ satisfies $\DDM^2\phi = \lambda_1^2 \phi.$
But, it is known that \cite{ginou}
\begin{eqnarray}\label{witten}
\DDM^2\phi = \hat D^2\phi +\frac n2 dH\cdot \nu\cdot\phi + \frac{n^2 H^2}{4}\phi,
\end{eqnarray}
where $\hat D$ is the Dirac-Witten operator given by $\hat D = \sum_{j=1}^{n} e_j\cdot\nabla^{\cZ}_{e_j}$. 
Hence, using that $\hat D\psi = -n\alpha\psi$ and (\ref{witten}), we get
$$\lambda_1^2 \phi = n^2\alpha^2 \phi +\frac n2 dH\cdot\nu\cdot\phi + \frac{n^2}{4} H^2 \phi.$$
Considering the real part of the scalar product of the last equality by $\phi$ implies that
$\lambda_1^2 = n^2\alpha^2 + \frac{n^2 H^2}{4}$. Hence, $H$ is constant.
\begin{example}
Simply connected complete Riemannian $\Spinc$ manifolds carrying parallel spinors were described in Examples \ref{Ex1}. 
The only $\Spinc$ structures on an irreducible K\"{a}hler not Ricci-flat manifold $\cZ$  which carry  
parallel spinors are the canonical and the anti-canonical one. In both cases, $\mu (\cZ, 0) = 1$ \cite{19}. Hence,
 Inequality (\ref{upperr}) holds for the first eigenvalue of the Dirac operator $\widetilde D$ defined on any
 compact Riemannian hypersurface endowed with the restricted $\Spinc$ structure. The complex projective space $\CC P^m$ 
or the complex hyperbolic space $\CC H^m$
with the Fubini-Study metric are examples of irreducible K\"{a}hler not Ricci-flat manifolds. 
From Examples \ref{Ex1}, the equality case in (\ref{lower}) 
is achieved for the first positive eigenvalue of the Dirac operator $\widetilde D$ defined
on Berger spheres embedded into $\CC P^2$. Hence, Inequality (\ref{upperr}) is also an equality in this case. 
Also, for rotational constant mean curvature $H$ spheres embedded into $\SSS^2 \times \RR$, 
Inequality (\ref{upperr}) is an equality because in this case, Inequality (\ref{lower}) is an equality. 
The only complete simply connected $\Spinc$ manifolds admitting real Killing spinors other than the $\Spin$ 
manifolds are the 
non-Einstein Sasakian manifolds endowed  with their canonical or anti-canonical 
$\Spinc$ structure \cite{19}. The manifolds $\mathbb{E}^* (\kappa, \tau)$ are examples of $\Spinc$ 
manifolds carrying a Killing spinor $\psi$ 
of Killing constant $\frac {\tau}{2}$. 
\end{example}
\section{Spin$^c$ structures on $\mathbb{E}^* (\kappa, \tau)$ and applications}\label{section4}
In this section, we make use of the existence of a $\Spinc$ Killing spinor to immerse $\mathbb{E}^* (\kappa, \tau)$ into complex space forms, 
to calculate some eigenvalues of the Dirac operator on Berger spheres and to prove the non-existence of totally umbilic surfaces in $\mathbb{E}^* (\kappa, \tau)$.
\subsection{Isometric immersions of $\mathbb{E}^* (\kappa, \tau)$ into complex space forms}
From the existence of an isometric embedding of $\mathbb{E}^* (\kappa, \tau)$ into $\M^4(\frac{\kappa}{4}- \tau^2)$, 
we reprove that the only $\Spinc$ structures on $\mathbb{E}^*(\kappa, \tau)$ carrying a Killing spinor 
are the canonical and the anti-canonical one. Conversely, 
the existence of a $\Spinc$ Killing spinor allows to immerse $\mathbb{E}^* (\kappa, \tau)$ in $\M^4(\frac{\kappa}{4}- \tau^2)$. More generally,
 we give necessary and sufficient geometric conditions to immerse any $3$-dimensional Sasaki manifold into $\M^2(c)$ for 
some $c\in\RR^*$.
\begin{prop}\label{im-inv}
The only $\Spinc$ structures on $\mathbb{E}^* (\kappa, \tau)$ carrying a real Killing spinor  are the canonical and the
 anti-canonical one. Moreover, the Killing constant is given by $\frac{\tau}{2}$.
\end{prop}
{\bf Proof.} It is known that there exists an isometric embedding of $\mathbb{E}^* (\kappa, \tau)$ into $\M^4(\frac{\kappa}{4}-\tau^2)$ of constant mean curvature $H = \frac{\kappa-16\tau^{2}}{12\tau}$ \cite{CMC}. Moreover, the second fundamental form is given by 
$$II(X) = -\tau X - \frac{4\tau^2-\kappa}{\tau} g_{\M^4}(X, \xi) \xi,$$
for every $X\in \Gamma(T \mathbb{E}^* (\kappa, \tau))$. Here, we recall that the normal vector of the immersion is given by $\nu :=J\xi$
 and $\{e_1, e_2, \xi, \nu = J\xi\}$ is a local orthonormal basis tangent to $\M^4$ where  $\{e_1, e_2, \xi\}$ 
is the canonical frame of $\mathbb{E}^* (\kappa, \tau)$. We denote by $\eta$ the real $1$-form associated with $\xi$, i.e., $\eta(X) = g(X, \xi)$ 
for any $X \in \Gamma(T\mathbb{E}^* (\kappa, \tau))$. The restriction of the canonical $\Spinc$ structure on $\M^4(\frac{\kappa}{4}-\tau^2)$
 induces a $\Spinc$ structure  on $\mathbb{E}^* (\kappa, \tau)$ and by the $\Spinc$  Gauss formula (\ref{spingauss}), the restriction of the 
parallel spinor on  $\M^4(\frac{\kappa}{4}-\tau^2)$ induces a 
spinor field $\phi$ on $\mathbb{E}^* (\kappa, \tau)$ satisfying, for all $X\in \Gamma(T\mathbb{E}^* (\kappa, \tau))$,
$$\nabla_X\phi = \frac{\tau}{2} X\bullet\phi  + \frac{4\tau^2-\kappa}{8\tau} \eta(X) \xi\bullet\phi.$$ 
Moreover, the spinor field $\phi$ satisfies $\xi\bullet\phi = -i\phi$ \cite[Theorem 3]{NR} and the curvature of the auxiliary 
line bundle $L$ associated with the induced $\Spinc$ structure is given by \cite[Theorem 3]{NR}
\begin{eqnarray}\label{omm}
 i\Omega(e_1, e_2) =-6i(\frac{\kappa}{4} -\tau^2), \ \ \ \ \text{and}\ \ \  \ i\Omega(e_i, e_j)=0,
\end{eqnarray}
elsewhere in the basis $\{e_1, e_2, \xi\}$. We deduce that, for all $X\in \Gamma(T\mathbb{E}^* (\kappa, \tau))$,
$$\nabla_X\phi = \frac{\tau}{2} X\bullet\phi - i\frac{4\tau^2-\kappa}{8\tau} g(X, \xi)\phi.$$
The connection $A$ on  the $\cercle^1$-principal fiber bundle $\cercle^1 (\mathbb{E}^* (\kappa, \tau))$ associated with the induced
 $\Spinc$ structure is the restriction to  $\mathbb{E}^* (\kappa, \tau)$ of the connection on the $\cercle^1$-principal fiber bundle $\cercle^1 \M^4$  associated 
with the canonical $\Spinc$ structure 
on $\M^4(\frac{\kappa}{4}-\tau^2)$, i.e.,  the connection $A$ on $\cercle^1 (\mathbb{E}^* (\kappa, \tau))$ is the restriction to $\mathbb{E}^* (\kappa, \tau)$ of the connection on $\cercle^1 (\M^4(\frac{\kappa}{4}-\tau^2))$ defined by 
the Levi-Civita connection. 
Let $\alpha$ be the real $1$-form on $\mathbb{E}^* (\kappa, \tau)$ defined by 
$$\alpha(X) = \frac{4\tau^2-\kappa}{4\tau} g(X, \xi),$$
for any $X\in \Gamma (T\mathbb{E}^* (\kappa, \tau))$. We endow the $\SSS^1$-principal fiber bundle $\SSS^1 (\mathbb{E}^* (\kappa, \tau))$ with the connection $A^{'} = A +i\alpha$. From (\ref{nnnabla}), 
there exists on $\Sigma \mathbb{E}^* (\kappa, \tau)$ a covariant derivative $\nabla^{'}$ such that
\begin{eqnarray*}
\nabla_X^{'}\phi = \nabla_X \phi + \frac{i}{2} \alpha(X)\phi = \frac{\tau}{2}  X\bullet\phi,
\end{eqnarray*}
for all $X\in \Gamma(T\mathbb{E}^* (\kappa, \tau))$.
Hence, we obtain a $\Spinc$ structure on $\mathbb{E}^* (\kappa, \tau)$ carrying a Killing spinor field and whose $\SSS^1$-principal 
fiber bundle $\SSS^1 (\mathbb{E}^* (\kappa, \tau))$ has a connection given by $A^{'}$. Now, we should prove that this $\Spinc$ structure 
is the canonical one. First, we calculate the curvature $i\Omega^{'} = i\Omega + id\alpha$ of $A^{'}$. 
It is easy to check that $\xi \lrcorner d\alpha = 0$ and 
$d\alpha (e_1, e_2) = -\frac{4\tau^2 - \kappa}{2}$. 
Hence, using (\ref{omm}), we get
$$\Omega^{'}(e_1, e_2) = - (\kappa-4\tau^2)\ \ \text{and}\ \ \xi \lrcorner \Omega^{'} = 0.$$
The curvature form $i\Omega^{'}$ is the same as the curvature form associated with the connection on the
 auxiliary line bundle of the canonical $\Spinc$ structure on $\mathbb{E}^* (\kappa, \tau)$. Since $\mathbb{E}^* (\kappa, \tau)$ is a simply
 connected manifold, we deduce that the $\SSS^1$-principal fiber bundle $\SSS^1 (\mathbb{E}^* (\kappa, \tau))$ endowed with 
the connection $A^{'}$ is the auxiliary  line bundle of the canonical $\Spinc$ structure on $\mathbb{E}^* (\kappa, \tau)$. 
Hence, we have on $\mathbb{E}^* (\kappa, \tau)$ two $\Spinc$ structures  with the same auxiliary line bundle 
(the canonical one and the one obtained by restriction of the canonical one on $\M^4$). 
But, on a Riemannian manifold $M$,  $\Spinc$ structures having the same auxiliary line bundle are parametrized 
by $H^1(M, \ZZ_2)$ \cite{mon}, which is trivial in our case since $\mathbb{E}^* (\kappa, \tau)$ is simply connected. 
To get the anti-canonical $\Spinc$ structure on $\mathbb{E}^* (\kappa, \tau)$, we restrict the anti-canonical $\Spinc$ structure on $\M^4$. 
In this case, $\xi\bullet\phi = i\phi$, $\Omega(e_1, e_2) = 6 (\frac{\kappa}{4} -\tau^2)$, $\xi\lrcorner\Omega =0$ and we 
choose the real $1$-form $\alpha$ to be $\alpha(X) = -\frac{4\tau^2-\kappa}{4\tau} g(X, \xi).$ 
\\\\
Next, we want to prove the converse. Indeed, we have
\begin{prop}\label{im}
 The manifolds $\mathbb{E}^*(\kappa, \tau)$ are isometrically immersed into $\M^4(c)$ for some $c$. Moreover, $\mathbb{E}^* (\kappa, \tau)$ are of constant mean curvature and $\eta$-umbilic.
\end{prop}
{\bf Proof.}  We recall that the $3$-dimensional homogeneous manifolds $\mathbb{E}^* (\kappa, \tau)$ have a $\Spinc$ structure 
(the canonical $\Spinc$ structure) carrying a Killing spinor field $\varphi$ of Killing constant 
$\frac{\tau}{2}$. Moreover $\xi\bullet\varphi = -i\varphi$ and 
\begin{eqnarray}
\Omega (e_1, e_2) = -(\kappa - 4 \tau^2) \ \ \text{and} \ \ \ \Omega(e_i, e_j) = 0,
\end{eqnarray}
in the basis $\{e_1, e_2, e_3 = \xi\}$. We denote by $A$ the connection on the auxiliary line bundle 
defining the canonical $\Spinc$ structure. 
Let $\alpha$ be the real $1$-form on $\mathbb{E}^* (\kappa, \tau))$ defined by  
$\alpha(X) = -\frac{4\tau^2-\kappa}{4\tau} g(X, \xi),$
for any $X\in \Gamma (T\mathbb{E}^* (\kappa, \tau))$. We endow the $\SSS^1$-principal fiber bundle $\SSS^1 (\mathbb{E}^* (\kappa, \tau))$ with the connection $A^{'} = A + i\alpha$. 
Then, there exists on $\Sigma \mathbb{E}^* (\kappa, \tau)$ a covariant derivative $\nabla^{'}$ such that
\begin{eqnarray}\label{inv}
\nabla_X^{'}\phi &=& \frac{\tau}{2}  X\bullet\phi + \frac{i}{2} \alpha(X)\phi \nonumber\\
&=& \frac{\tau}{2}  X\bullet\phi  + \frac{4\tau^2-\kappa}{8\tau} \eta(X)\xi\bullet\phi,
\end{eqnarray}
for all $X\in  \Gamma(T\mathbb{E}^* (\kappa, \tau))$.
Hence, we obtain a $\Spinc$ structure on $\mathbb{E}^* (\kappa, \tau)$ carrying a 
 spinor field $\phi$ satisfying (\ref{inv}) and whose $\SSS^1$-principal fiber bundle $\SSS^1 (\mathbb{E}^* (\kappa, \tau))$ has a connection
 given by $A^{'}$. 
We calculate the curvature $i\Omega^{'} = i\Omega + id\alpha$ of $A^{'}$. It is easy to check that $\xi \lrcorner d\alpha = 0$ and 
$d\alpha (e_1, e_2) = \frac{4\tau^2 - \kappa}{2}$. Hence,
$$\Omega^{'}(e_1, e_2) = - 6(\frac{\kappa}{4}-\tau^2)\ \ \text{and}\ \ \xi \lrcorner \Omega^{'} = 0.$$
Since $\mathbb{E}^* (\kappa, \tau)$ are Sasakian, by \cite[Theorem $4$]{NR}, we get an isometric immersion of $\mathbb{E}^* (\kappa, \tau)$ into 
$\M^4(c)$ for $c =\frac{\kappa}{4}-\tau^2$. Moreover, $\mathbb{E}^* (\kappa, \tau)$ are of constant 
mean curvature  and $\eta$-umbilic (see \cite{NR}).\\\\
More general, we have:
\begin{thm}\label{sasaki-immersion}
 Every  simply connected  non-Einstein  $3$-dimensional Sasaki manifold $M^3$ of constant scalar curvature can be immersed into
 $\M^4(c)$ for some $c \in \RR^*$. Moreover, $M$ is $\eta$-umbilic and of constant mean curvature.
\end{thm}
{\bf Proof. } We recall that a Sasaki structure  on a $3$-dimensional manifold $M^3$ is given by a Killing vector field $\xi$ of unit length such that the tensors 
$\Chi := \nabla\xi$ and $\eta := g(\xi, \cdot)$ are related by 
$$\Chi^2 = -\id + \eta\otimes \xi.$$
We know that a  non-Einstein Sasaki  manifold has a $\Spinc$ structure 
carrying a Killing 
spinor field $\phi$ of Killing constant $\beta$. By rescaling the metric, we can assume that $\beta= -\frac 12$. Moreover, 
the Killing vector field $\xi$ defining the Sasaki structure satisfies $\xi\bullet\phi = -i\phi$ (see \cite{19}). The Ricci tensor on $M$ is given by
$$\mathrm{Ric}(e_j) = \frac{S-2}{2}e_j, j=1,2\ \ \ \ \ \text{and}\ \ \ \ \mathrm{Ric}(\xi) = 2\xi,$$
where $S$ denotes the scalar curvature of $M$ and $\{e_1, e_2, \xi\}$  a local orthonormal frame of $M$.
Because we assumed that $M$ is non-Einstein, we have $S \neq 6$ and hence we can find $c\in \RR^*$ such that $S = 8c+6$. 
The Ricci identity (\ref{rici}) in $X = \xi$ gives that $\xi\lrcorner\Omega = 0$ and by the Schr\"{o}dinger-Lichnerowicz formula, it follows that
$\Omega(e_1, e_2) = \frac{6-S}{2}$. Let $\alpha$ be the real $1$-form on $M$ defined by  
$\alpha(X) = -c g(X, \xi),$
for any $X\in \Gamma (TM)$. We endow the $\SSS^1$-principal fiber bundle $\SSS^1M$ with the connection $A^{'} = A +i\alpha$, where $A$ denotes the
 connection on $\SSS^1 M$ whose curvature form is given by $i\Omega$.
From (\ref{nnnabla}), there exists on $\Sigma M$ a covariant derivative $\nabla^{'}$ such that
\begin{eqnarray*}
\nabla^{'}_X\phi =-\frac{1}{2}  X\bullet\phi  - \frac{i}{2} c g(X, \xi)\phi.
\end{eqnarray*}
Now, we calculate the curvature $i\Omega^{'} = i\Omega + id\alpha$ of $A^{'}$. It is easy to check that $\xi \lrcorner d\alpha = 0$ and 
$d\alpha (e_1, e_2) = -2c$. Hence,
$$\xi\lrcorner\Omega^{'} = 0, \ \ \ \ \ \Omega^{'}(e_1, e_2) = \frac{6-S}{2}-2c = -6c.$$ 
By \cite[Theorem $4$]{NR}, $M$ is immersed into $\M^4(c)$. Additionally, $M$ is $\eta$-umbilic and of constant mean curvature. 
\subsection{Totally umbilic surfaces in $\mathbb{E}^* (\kappa, \tau)$}
By restriction of the Killing spinor of Killing constant $\frac{\tau}{2}$ on $\mathbb{E}^* (\kappa, \tau)$
to a surface $M^2$, the authors characterized isometric immersions into $\Ekt$ by the existence of a $\Spinc$ 
structure carrying a special spinor field  \cite{NR}. More precisely, consider
 $(M^2, g)$  a Riemannian surface. We denote by  $E$ a field of symmetric endomorphisms of $TM$, 
with trace equal to $2H$. The vertical vector field can be written as  $\xi=dF(T)+f\nu$, where $\nu$ is the unit normal vector to
 the surface, $f$ is a real function on $M$ and $T$ the tangential part of $\xi$. The 
isometric immersion of $(M^2,g)$ into $\Ekt$ with shape operator $E$, mean curvature $H$ is characterized by a $\Spinc$ structure on $M$ carrying a non-trivial spinor field $\phi$ satisfying,
 for all $X\in \Gamma(TM)$,
$$\nabla_X\phi=-\frac{1}{2}EX\bullet\phi+i\frac{\tau}{2}X\bullet\overline{\phi}.$$
Moreover, the auxiliary bundle has a connection of curvature given, in any local orthonormal 
frame $\{t_1, t_2\}$, by   
$i\Omega(t_1, t_2)= -i(\kappa-4\tau^2)f =-i(\kappa-4\tau^2) \frac{<\phi, \overline \phi>}{\vert\phi\vert^2} $. The vector $T$ is given by 
$$g(T, t_1) = <it_2\bullet\phi, \frac{\phi}{\vert\phi\vert^2}> \ \ \ \text{and}\ \ \ \ g(T, t_2) = -<it_1\bullet\phi, \frac{\phi}{\vert\phi\vert^2}>.$$
Here and also by restriction of the Killing spinor, we gave an elementary $\Spinc$ proof of the following result proved by  R. Souam and E. Toubiana in \cite{ST}.
\begin{thm}\label{tum}
There are no totally umbilic surfaces in $\mathbb{E}^*(\kappa, \tau)$.
\end{thm}
{\bf Proof.} Assume that $M$ is a totally umbilical surface in $\mathbb{E}^* (\kappa, \tau)$, i.e. $E = H\ \id$. Then $d^\nabla E(e_1, e_2) = (\nabla_{t_1}E) t_2 - (\nabla_{t_2}E) t_1 = J(dH)$. The $\Spinc$ curvature $\mathcal{R}$ on the spinor  field $\phi$  is given by  \cite{NR}:
$$\mathcal R(t_1, t_2)\phi = -\frac 12 J(dH)\bullet\phi + i\frac{H^2}{2}\overline \phi + i\frac{\tau^2}{2}\overline \phi.$$
The $\Spinc$ Ricci identity (\ref{rici}) on the surface $M$ implies
\begin{eqnarray*}
t_1\bullet \mathcal R(t_1, t_2)\phi = \frac 12 \mathrm{Ric}(t_2)\bullet\phi -\frac i2 (t_2\lrcorner \Omega)\bullet\phi
\end{eqnarray*}
Hence,
\begin{eqnarray*}
-\frac 12 t_1\bullet J(dH)\bullet\varphi + \frac i2 H^2 t_1\bullet\overline \phi + \frac i2 \tau^2 t_1\bullet \overline \phi  = \frac 12 \mathrm{Ric}(t_2)\bullet\phi +\frac i2  \Omega(t_1, t_2) t_1\bullet\phi
\end{eqnarray*}
Consider the real part of the scalar product of the last identity by $\phi$, we get
$$g(t_1, J(dH)) = \Omega(t_1, t_2) <i t_1\bullet\phi, \frac{\phi}{\vert\phi\vert^2}> =  -\Omega(t_1, t_2) g(T, t_2).$$
Finally, $-g(t_2, dH) = (\kappa-4\tau^2) f g(T, t_2)$. The same holds for $t_1$. Then, 
$$dH = -(\kappa-4\tau^2) f T,$$ 
which gives the contradiction. The last identity  is the same obtained by  R. Souam and E. Toubiana in \cite{ST}.
\subsection{Spectrum of the Spin$^c$ Dirac operator on Berger spheres}
In this subsection, we apply a method of C.  B\"ar \cite{bar96, ginoulivre} to find 
explicitly some eigenvalues of the $\Spinc$ Dirac operator on Berger spheres, i.e., on $\mathbb{E}^* (\kappa, \tau)$ with $\kappa > 0$.
\begin{lemma}
Let $(M^n, g)$ be a Riemannian $\Spinc$ manifold carrying a Killing spinor $\phi$ of Killing number  $\alpha\in \RR^*$. Then, 
$(\lambda_k(\bigtriangleup)  + (\frac{n-1}{2})^2)_{k\in \NN}$ are some eigenvalues of $(D+ \frac{\alpha}{2} \id)^2$. Here $\lambda_k(\bigtriangleup)$, $k=0,1,...$  denote the eigenvalue of the Laplacian $\bigtriangleup$.
\label{spe}
\end{lemma}
{\bf Proof:} We have $D\phi = -\frac{n\alpha}{2} \phi$. For every $f\in C^\infty (M, \RR)$, we can easily check that
\begin{eqnarray*}
D^2(f\phi) = (\frac{n^2}{4} -\frac{n}{2})f\phi - \alpha D(f\phi)+ (\bigtriangleup f)\phi,
\end{eqnarray*}
Hence, $(D+ \frac{\alpha}{2} \id)^2 (f\phi) = ( \bigtriangleup f + (\frac{n-1}{2})^2f)\phi.$ Now, if $\{f_k\}_{k\in \NN}$ denotes a $L^2$-orthonormal basis of eigenfunctions of $\bigtriangleup$ of $M$, then for every $k\in \NN$, we get
$$(D+ \frac{\alpha}{2} \id)^2 (f_k\phi) = ( \lambda_k(\bigtriangleup)  + (\frac{n-1}{2})^2)f_k\phi,$$
where $\lambda_k(\bigtriangleup)$ is the eigenvalue 
of $\bigtriangleup$ whose eigenfunction is $f_k$. So, $(\lambda_k(\bigtriangleup)  
+ (\frac{n-1}{2})^2)_{k\in \NN}$ are some eigenvalues of $(D+ \frac{\alpha}{2} \id)^2$.
\\\\
{\bf Spectrum of Berger spheres endowed with the canonical Spin$^c$ structure.}
We consider Berger spheres with Berger metrics $g_{\kappa, \tau}$, $\kappa > 0$ and $\tau \neq 0$ defined by 
$$g_{(\kappa, \tau)} (X, Y) =  \frac{\kappa}{4}\Big(g(X, Y) + (\frac{4\tau^2}{\kappa}-1) g(X, \xi)g(Y, \xi)\Big),$$
where $g$ is the standard metric on $\SSS^3$ of constant curvature $1$.
For simplicity, we can assume that $\kappa =4$ ($\tau\neq \pm 1$). For any function $f$, the Laplacian $\bigtriangleup_{4, \tau}$ with respect to $g_{4, \tau}$ is related to the Laplacian $\bigtriangleup$ with respect to $g$ by \cite{tanno}
$$\bigtriangleup_{4, \tau} f = \bigtriangleup f -(1-\tau^{-2})\xi(\xi(f)).$$
It is known that each eigenfunction $f_k$ of $\bigtriangleup$ corresponding to $\lambda_k (\bigtriangleup) = k(2+k)$ ($k \in \mathbb{N}$) is also an eigenfunction of $\bigtriangleup_{4, \tau}$ \cite{tanno} corresponding to 
$$\lambda_k(\bigtriangleup) -(1-\tau^{-2})(k-2p)^2,\ \ \ \ 0\leq p\leq [\frac{k}{2}].$$
Moreover, each eigenvalue of $\bigtriangleup_{4,\tau}$ takes the above form. We recall that the eigenspace 
of $\bigtriangleup$ corresponding to $\lambda_k(\bigtriangleup)$  is the restriction to the sphere $\SSS^3$ of the set of  harmonic homogeneous polynomial on $\RR^4$ of degree $k$. When we consider Berger spheres endowed with the canonical $\Spinc$ structure, we get by Lemma \ref{spe}
$$\Big(D+\frac{\tau}{2}\id\Big)^{2} (f_k \phi)= \Big[2+k(2+k)- (1-\tau^{-2}) (k-2p)\Big] f_k\phi,$$
where $\phi$ is the Killing spinor field of Killing constant $\frac{\tau}{2}$. Hence,
$$\mu_{k,p} = -\frac{\tau}{2} \pm \sqrt{2+k(2+k)-(1-\tau^{-2})(k-2p)}$$
are some eigenvalues of the Dirac operator on Berger spheres with $-1 <\tau < 1$ and  endowed with the canonical $\Spinc$ structure.
\\\\
{\bf Spectrum of Berger spheres endowed with the Spin$^c$ structure induced from the canonical one on $\M^4(1-\tau^2)$.}
On Berger spheres, we have shown that the $\Spinc$ structure induced from the canonical one on $\M^4 (1- \tau^2)$ carries a spinor field $\phi$ satisfying 
\begin{eqnarray*}
\nabla_X\phi = \frac{\tau}{2} X\bullet\varphi - i\frac{\tau^2-1}{2\tau}g(X, \xi)\phi = \nabla^{'}_X \phi - i\frac{\tau^2-1}{2\tau}g(X, \xi)\phi
\end{eqnarray*}
Then, denoting by $D$ (resp. $D^{'}$) the Dirac operator associated with the restricted $\Spinc$ structure (resp. with the canonical $\Spinc$ structure), we get $D\phi = D^{'}\phi-\frac{\tau^2-1}{2\tau}\phi$. for any function $f$, we have
\begin{eqnarray*}
D(f\phi)= \grad f \cdot\phi + f D\phi = D^{'} (f\phi) - f D^{'}\phi + fD\phi =  D^{'} (f\phi) - (\frac{\tau^2-1}{2\tau})f\phi.
\end{eqnarray*}
Hence, we have $D(f_k\phi) = \Big(\mu_{k,p} -\frac{\tau^2-1}{2\tau}\Big) f_k\phi$ and $\mu_{k,p} -\frac{\tau^2-1}{2\tau}$ are some eigenvalues of the Dirac operator on Berger spheres endowed with the $\Spinc$ structure induced from the canonical one on $\M^4(1-\tau^2)$, $-1 < \tau <1$.
\\\\{\bf Acknowledgment.} Both authors would like to thank Sebasti\'an Montiel for pointing out this  correspondence between Killing spinors on $\mathbb E^*(\kappa, \tau)$ and isometric immersions into $\M^4(c)$.  Also, we would like to thank  Oussama Hijazi for helpful discussions and relevant remarks.

\end{document}